\documentclass[10pt,letterpaper,english]{article}
\usepackage[T1]{fontenc}
\usepackage[latin9]{inputenc}
\usepackage{xcolor}
\usepackage{babel}
\usepackage{units}
\usepackage{amsmath}
\usepackage{amssymb}
\usepackage[all]{xy}
\PassOptionsToPackage{normalem}{ulem}
\usepackage{ulem}
\usepackage[unicode=true,pdfusetitle,
 bookmarks=true,bookmarksnumbered=false,bookmarksopen=false,
 breaklinks=false,pdfborder={0 0 1},backref=section,colorlinks=false]
 {hyperref}
\usepackage{breakurl}

\makeatletter

\providecolor{lyxadded}{rgb}{0,0,1}
\providecolor{lyxdeleted}{rgb}{1,0,0}

\@ifundefined{date}{}{\date{}}
\usepackage{babel}
\usepackage{cite}

\usepackage{color}

\usepackage{lineno}

\usepackage[labelfont=bf,labelsep=period,justification=raggedright]{caption}

\renewcommand{\@biblabel}[1]{\quad#1.}

\usepackage{placeins}

\AtBeginDocument{
  
}

\makeatother

\begin{document}
\begin{flushleft}
\textbf{\large{Embedding Quantum into Classical: }}\\
\textbf{\large{Contextualization vs Conditionalization}}{\large{}}\\
{\large{ Ehtibar N. Dzhafarov$^{1\ast}$, Janne V. Kujala$^{2}$}}\\
{\large{ }}\textbf{\large{1}}{\large{ Department of Psychological
Sciences, Purdue University, West Lafayette, Indiana, USA}}\\
{\large{ }}\textbf{\large{2}}{\large{ Department of Mathematical Information
Technology, University of Jyväskylä, Jyväskylä, Finland}}\\
{\large{ $\ast$ E-mail: ehtibar@purdue.edu}}
\par\end{flushleft}{\large \par}

\section*{Abstract}

We compare two approaches to embedding joint distributions of random
variables recorded under different conditions (such as spins of entangled
particles for different settings) into the framework of classical,
Kolmogorovian probability theory. In the contextualization approach
each random variable is ``automatically'' labeled by all conditions
under which it is recorded, and the random variables across a set
of mutually exclusive conditions are probabilistically coupled (imposed
a joint distribution upon). Analysis of all possible probabilistic
couplings for a given set of random variables allows one to characterize
various relations between their separate distributions (such as Bell-type
inequalities or quantum-mechanical constraints). In the conditionalization
approach one considers the conditions under which the random variables
are recorded as if they were values of another random variable, so
that the observed distributions are interpreted as conditional ones.
This approach is uninformative with respect to relations between the
distributions observed under different conditions, because any set
of such distributions is compatible with any distribution assigned
to the conditions. \thispagestyle{empty} \markright{Contextualization
vs Conditionalization}

\section*{Introduction}

\subsection*{Joint Distributions and Stochastic Unrelatedness}

Many scientific problems, from psychology to quantum mechanics, can
be presented in terms of \emph{random outputs} of some system recorded
under various \emph{conditions}. According to the principle of \emph{Contextuality-by-Default}
{[}\ref{Dzhafarov,-E.N.,-&2013}-\ref{Dzhafarov,-E.N.,-&New}{]},
when applying Kolmogorov's probability theory (KPT) to such a problem,
random variables recorded under different, mutually incompatible conditions
should be viewed as \emph{stochastically unrelated} to each other,
i.e., possessing no joint distribution. They can always be ``sewn
together'' as part of their theoretical analysis, but joint distributions
are then \emph{imposed} on them rather than derived from their identities.
In this paper we discuss two possible approaches to the foundational
issue of ``sewing together'' stochastically unrelated random variables.
We call these approaches \emph{contextualization} and \emph{conditionalization}.
The former takes the Contextuality-by-Default principle as its departure
point and is, in a sense, its straightforward extension; in the latter,
Contextuality-by-Default is obtained as a byproduct. 

To understand why the Contextuality-by-Default principle is associated
with either of these two approaches, one should first of all abandon
the naive notion that in KPT any two random variables have a joint
distribution uniquely determined by their definitions. A random variable
is a measurable function on a probability space, and the notion of
a single probability space for all possible random variables (or,
equivalently, the notion of a single random variable of which all
other random variables are functions)%
\footnote{In this discussion we impose no restrictions on the domain and codomain
probability spaces. A random variable therefore is understood in the
broadest possible way, including random vectors, random functions,
random sets, etc. We will avoid, however, the use of general measure-theoretic
formalism.%
} is untenable. It contradicts the commonly used KPT constructions. 

One of them is, given any set, to construct a random variable whose
range of possible values coincides with this set. A probability space
on which all such random variables were defined would have to include
a set of cardinality exceeding that of all possible sets, an impossibility. 

Another commonly used construction is, given a random variable, to
introduce another random variable that has a given distribution and
is stochastically independent of the former. The use of this construction
contradicts even the notion of a jointly distributed set of all variables
with a particular distribution {[}\ref{Dzhafarov,-E.N.,-&LNCS13}{]},
say, the set $\mathrm{Norm}\left(0,1\right)$ of all standard-normally
distributed random variables. Indeed, if all random variables in $\mathrm{Norm}\left(0,1\right)$
were jointly distributed, they would all be presentable as functions
of some random variable $N$, the identity function on the probability
space on which the random variables in $\mathrm{Norm}\left(0,1\right)$
are defined. Choose now a standard-normally distributed random variable
$X$ so that it is independent of $N$. Then it is also independent
of any $Y\in\mathrm{Norm}\left(0,1\right)$. Since $X$ cannot be
independent of itself, $X$ cannot belong to $\mathrm{Norm}\left(0,1\right)$.
At the same time, $X$ must belong to $\mathrm{Norm}\left(0,1\right)$
due to its distribution. 

Short of imposing on KPT artificial constraints (such as an upper
limit on cardinality of the random variables' ranges), these and similar
contradictions can only be dissolved by allowing for stochastically
unrelated random variables defined on different probability spaces
(see Ref. {[}\ref{Dzhafarov-EN,-Kujala_Chapter}{]} for how this can
be built into the basic set-up of probability theory). The principle
of Contextuality-by-Default eliminates guesswork from deciding which
random variables are and which are not jointly distributed. Irrespective
of how one defines a system with random outputs and identifies the
conditions under which these outputs are recorded, the outputs are
jointly distributed if they are recorded under one and the same set
of conditions; if they are recorded under different, mutually exclusive
conditions, they are stochastically unrelated.

\subsection*{Two Approaches }

Contextualization and conditionalization differ in how they ``sew
together'' stochastically unrelated random variables. To demonstrate
these differences on a simple example, let $X$ and $Y$ be random
variables with $+1/-1$ values, so that their distributions are determined
by $\Pr\left[X=1\right]$ and $\Pr\left[Y=1\right]$, respectively.
Let $X$ and $Y$ be recorded under mutually exclusive conditions. 

In contextualization (the approach we proposed in Refs. {[}\ref{Dzhafarov,-E.N.,-&2013}-\ref{Dzhafarov,-E.N.,-&New}{]}),
one first invokes the Contextuality-by-Default principle to treat
$X$ and $Y$ as stochastically unrelated random variables. A ``sewing
together'' of $X$ and $Y$ consists in \emph{probabilistically coupling}
them {[}\ref{Thorisson,-H.:-Coupling,}{]}, i.e., presenting them
as functions of a single random variable. Put differently (but equivalently),
we create a random variable (vector) $Z=\left(X',Y'\right)$ such
that $X'$ is distributed as $X$ and $Y'$ is distributed as $Y$.
The variables $X'$and $Y'$ are jointly distributed (otherwise $Z$
would not be called a random variable, or a random vector), but this
distribution is not unique. Thus, $X$ and $Y$ can always be coupled
as stochastically independent random variables, so that 
\begin{equation}
\Pr\left[X'=1,Y'=1\right]=\Pr\left[X'=1\right]\times\Pr\left[Y'=1\right].
\end{equation}
They can also be coupled as identical random variables, 
\begin{equation}
\Pr\left[X'=Y'\right]=1,
\end{equation}
but only if $X$ and $Y$ are distributed identically,
\begin{equation}
\Pr\left[X=1\right]=\Pr\left[Y=1\right].
\end{equation}
There can, in fact, be an infinity of couplings, constrained only
by 
\begin{equation}
\begin{array}{c}
\Pr\left[X'=1,Y'=1\right]+\Pr\left[X'=1,Y'=-1\right]=\Pr\left[X=1\right],\\
\Pr\left[X'=1,Y'=1\right]+\Pr\left[X'=-1,Y'=1\right]=\Pr\left[Y=1\right].
\end{array}
\end{equation}

In conditionalization, one creates a random variable $C$ with two
possible values corresponding to the two sets of conditions under
which one records $X$ and $Y$, respectively. Then one defines a
random variable $U=\left(C,V\right)$, such that the conditional distribution
of $V$ given $C=1$ is the same as the distribution of $X$, and
the conditional distribution of $V$ given $C=2$ is the same as the
distribution of $Y$. In other words,
\begin{equation}
\begin{array}{c}
\Pr\left[V=1\,\vert\, C=1\right]=\Pr\left[X=1\right],\\
\Pr\left[V=1\,\vert\, C=2\right]=\Pr\left[Y=1\right].
\end{array}
\end{equation}
The principle of Contextuality-by-Default here does not have to be
invoked explicitly, but it is adhered to anyway: the random variable
$V$ is related to conditions under which it is recorded, and $V$
conditioned on $C=1$ clearly has no joint distribution with $V$
conditioned on $C=2$. 

Conditionalization can also be implemented in more complex constructions,
such as the one proposed in Ref. {[}\ref{enu:D.-Avis,-P.}{]}. In
our example, this construction amounts to replacing $V$ with two
random variables, $V_{1}$ and $V_{2}$, and ``coordinating'' their
possible values with the values of $C$. Thus, one can make $V_{1}$
and $V_{2}$ binary, $+1/-1$, and define the conditional distributions
by
\begin{equation}
\begin{array}{c}
\Pr\left[V_{1}=v,V_{2}=1\,\vert\, C=1\right]=\Pr\left[V_{1}=v,V_{2}=-1\,\vert\, C=1\right]=\frac{1}{2}\Pr\left[X=v\right],\\
\Pr\left[V_{1}=1,V_{2}=v\,\vert\, C=2\right]=\Pr\left[V_{1}=-1,V_{2}=v\,\vert\, C=2\right]=\frac{1}{2}\Pr\left[Y=v\right],
\end{array}
\end{equation}
where $v=1$ or $-1$. For $C=1$, as we see, the ``relevant'' output
is $V_{1}$, and the probabilities of its values $v$ are simply evenly
divided between the two possible values of the ``irrelevant'' output
$V_{2}$ (and for $C=2$, $V_{1}$ and $V_{2}$ exchange places). 

We argue in this paper that only contextualization serves as a useful
tool for classifying and characterizing different types of systems
involving\emph{ }random outputs that depend on\emph{ }conditions (e.g.,
classical-mechanical vs quantum-mechanical systems). Conditionalization,
both in its simplest and modified versions, is always applicable but
uninformative.

\subsection*{Quantum Entanglement}

Our analysis pertains to any input-output relations, as considered
in Refs. {[}\ref{Dzhafarov,-E.N.,-&2013}-\ref{Dzhafarov,-E.N.,-&Advances},\ref{enu:Dzhafarov-EN,-Kujala2010}-\ref{enu:Dzhafarov-EN,-KujalaAMS}{]}.
The relations can be physical, biological, behavioral, social, etc.
For the sake of mathematical transparency, however, we confine our
consideration to the canonical quantum-mechanical paradigm {[}\ref{enu:Bohm-D,-Aharonov}{]}
involving two entangled particles, ``Alice's'' and ``Bob's.''
Alice measures the spin of her particle in one of two directions,
$\alpha_{1}$ or $\alpha_{2}$ (values of the first input), and Bob
measures the spin of his particle in one of two directions, $\beta_{1}$
or $\beta_{2}$ (values of the second input). Each pair of measurements
is therefore characterized by one of four possible combinations of\emph{
}input values $\left(\alpha_{i},\beta_{j}\right)$, and it is these
combinations that form the four \emph{conditions} in this example.
The spins recorded in each trial are realizations of random variables
(outputs) $A$ and $B$, which, in the simplest case of spin-$\nicefrac{1}{2}$
particles, can attain two values each: ``up'' or ``down'' (encoded
by $+1$ and $-1$, respectively). 

Aside from simplicity, another good reason for using this example
is that it relates to the problem of great interest in the foundation
of physics: in what way and to what an extent one can embed joint
probabilities of spins in entangled particles into the framework of
KPT? It may seem that this question was answered by John Bell in his
classical papers {[}\ref{Bell,-J.:-On},\ref{enu:Bell-J-(1966)}{]},
and that the answer was: KPT is not compatible with the joint distributions
of spins in entangled particles. However, in Bell's analysis and its
subsequent elaborations {[}\ref{Clauser,-J.F.,-Horne,},\ref{enu:Fine,-A.:-Hidden}{]}
the use of KPT is constrained by an added assumption that has nothing
to do with KPT. Namely, the implicit assumption in these analyses
is that of ``noncontextuality'': 
\begin{quote}
a spin recorded in Alice's particle is a random variable uniquely
identified by the measurement setting (spatial axis) for which it
is recorded (and analogously for Bob's particle). 
\end{quote}
In other words, the spin recorded by Alice for settings $\alpha_{1}$
and $\alpha_{2}$ are different random variables $A_{1}$ and $A_{2}$,
but the identity of either of them does not depend on whether Bob's
setting is $\beta_{1}$ or $\beta_{2}$ (and analogously for Bob's
random variables $B_{1},B_{2}$ corresponding to $\beta_{1}$ and
$\beta_{2}$). For well-established reasons (discussed in detail below),
this makes a Kolmogorovian account of quantum entanglement impossible.

However, according to the Contextuality-by-Default principle, if one
applies it to the Alice-Bob paradigm, 
\begin{quote}
any two random variables recorded under mutually exclusive conditions
are labeled by these conditions and considered stochastically unrelated\emph{.} 
\end{quote}
Alice's spin values recorded under the condition $\left(\alpha_{1},\beta_{1}\right)$
cannot \emph{co-occur} with the spin values recorded by her under
the condition $\left(\alpha_{1},\beta_{2}\right)$, even though $\alpha_{1}$
is the same in both conditions. Therefore the identity of the spin
she measures under $\left(\alpha_{1},\beta_{1}\right)$ should be
viewed as different from the identity of the spin she measures under
$\left(\alpha_{1},\beta_{2}\right)$. 

This leads one to the double-indexation of the spins, 
\begin{equation}
A{}_{11},A{}_{12},A{}_{21},A{}_{22},B{}_{11},B{}_{12},B{}_{21},B{}_{22},\label{eq:8}
\end{equation}
where $A_{ij}$ and $B_{ij}$ are the measurements by Alice and Bob,
respectively, recorded under the condition $\left(\alpha_{i},\beta_{j}\right)$,
$i,j\in\left\{ 1,2\right\} $. This vector of random variables cannot
be called a \emph{random vector} (or random variable, as we use the
term broadly), because its components are not jointly distributed.
Thus, $A{}_{11}$ and $A_{12}$, or $A{}_{11}$ and $B_{12}$, are
recorded under mutually exclusive conditions, so they do not have
jointly observed realizations. But the outputs $A_{11}$ and $B_{11}$,
being recorded under one and the same condition $\left(\alpha_{1},\beta_{1}\right)$,
are jointly distributed, i.e., the joint probabilities for different
combinations of co-occurring values of $A_{11}$ and $B_{11}$ are
well-defined. The situation is summarized in the following diagram:
\begin{equation}
\begin{array}{c}
\xymatrix{ & B_{12}\ar[dl] & B_{22}\ar[dr]_{joint}\\
A_{12}\ar[ur]_{joint} &  &  & A_{22}\ar[ul]\\
A_{11}\ar[dr]^{joint} &  &  & A_{21}\ar[dl]\\
 & B_{11}\ar[lu] & B_{21}\ar[ur]^{joint}
}
\end{array}\label{eq:all8}
\end{equation}

\subsection*{Contextuality and No-Signaling}

Why do we speak of ``contextuality'' and ``noncontextuality''?
The terms come from quantum mechanics (see, e.g., Refs. {[}\ref{enu:S.-Kochen,-F.}-\ref{enu:A.-Cabello:-Simple}{]}),
although it is not always clear that they are used in the same meaning
as in the present paper. In the Alice-Bob paradigm with two spin-$\nicefrac{1}{2}$
particles, the (marginal) distribution of Alice's measurement $A_{ij}$
does not depend on Bob's setting $\beta_{j}$, nor does the distribution
of Bob's measurement $B_{ij}$ depend on $\alpha_{i}$: 
\begin{equation}
\begin{array}{cc}
\Pr\left[A_{11}=1\right]=\Pr\left[A_{12}=1\right], & \Pr\left[A_{21}=1\right]=\Pr\left[A_{22}=1\right],\\
\Pr\left[B_{11}=1\right]=\Pr\left[B_{21}=1\right], & \Pr\left[B_{12}=1\right]=\Pr\left[B_{22}=1\right].
\end{array}\label{eq:marginal selectivity}
\end{equation}
This is known as the \emph{no-signaling} condition {[}\ref{Cereceda,-J.-Quantum}{]}:
Alice, by watching outcomes of her measurements, is not able to guess
Bob's settings, and vice versa. If the two particles are separated
by a space-like interval, violations of no-signaling would contravene
special relativity (and imply the ``spooky action at a distance,''
in Einstein's words). 

Nevertheless, in KPT, $A$ cannot be indexed by $\alpha_{i}$ alone,
nor can $B$ be indexed by $\beta_{j}$ alone. 

The logic forbidding single-indexation of the spins, $A_{1},A_{2},B_{1},B_{2}$,
is simple {[}\ref{Dzhafarov,-E.N.,-&New}{]}. Since, for any $i,j$,
the random variables $A_{i}$ and $B_{j}$ are jointly distributed,
they are defined on the same probability space. Applying this consideration
to $\left(A_{1},B_{1}\right)$, $\left(A_{1},B_{2}\right)$, and $\left(A_{2},B_{1}\right)$,
we are forced to accept that all four random variables, $A_{1},A_{2},B_{1},B_{2}$,
are defined on one and the same probability space. The existence of
this joint distribution, however, is known to be equivalent to Bell-type
inequalities (see below), known not to hold for entangled particles. 

Therefore, in perfect compliance with the Contextuality-by-Default
principle, we are forced to use the double indexation (\ref{eq:8}).
We can say that while $\beta_{j}$ does not influence $A_{ij}$ ``directly''
(which would be the case if $\beta_{j}$ could affect the distribution
of $A_{ij}$), it generally creates a ``context'' for $A_{ij}$.
The context makes $A_{i1}$ and $A_{i2}$ two different random variables
with one and the same distribution, rather than one and the same random
variable. (Analogous reasoning applies to $B_{ij}$ in relation to
$\alpha_{i}$.) 

It should not, of course, come as a surprise that different random
variables can be identically distributed. After all, it is perfectly
possible that the distributions of Alice's spins for $\alpha_{1}$
and $\alpha_{2}$ are identical too, and this would not imply that
they are one and the same random variable. Within the framework of
KPT, the difference between $A_{11}$ and $A_{12}$ is essentially
the same as the difference between $A_{11}$ and $A_{21}$: in both
cases we deal with stochastically unrelated random variables, the
only difference being that in the former pair, unlike in the latter
one, the no-signaling condition forces the two random variables to
be identically distributed. The notion of contextuality, however,
does require broadening of one's thinking about how one decides that
some empirical observations are and some are not realizations of one
and the same random variable, as understood in KPT {[}\ref{Dzhafarov,-E.N.,-&LNCS13},\ref{Dzhafarov,-E.N.,-&Advances}{]}.

\section*{Theory}

\subsection*{Contextualization and Couplings}

Contextualization is a straightforward extension of the Contextuality-by-Default
principle. The latter creates the eight random variables in (\ref{eq:8}),
and the contextualization approach consists in directly imposing a
joint distribution on them. This can, of course, be done in infinitely
many ways. Any random variable 
\begin{equation}
Y=\left(A'_{11},A'_{12},A'_{21},A'_{22},B'_{11},B'_{12},B'_{21},B'_{22}\right)\label{eq:coupling8}
\end{equation}
such that, for any $i=1,2$ and $j=1,2$, 
\begin{equation}
\left(A'_{ij},B'_{ij}\right)\textnormal{ is distributed as }\left(A_{ij},B_{ij}\right),\label{eq:coupling8 condition}
\end{equation}
is called a (probabilistic) \emph{coupling} for (\ref{eq:8}) {[}\ref{Thorisson,-H.:-Coupling,}{]}.
The fact that $Y$ in (\ref{eq:coupling8}) is referred to as a \emph{random
variable} (or random vector) implies that the components of $Y$ are
jointly distributed, i.e., there is a joint probability assigned to
each of the $2^{8}$ combinations of values for $A'_{11},A'_{12},\ldots,B'_{22}$.

The set of all possible couplings (\ref{eq:coupling8}) for (\ref{eq:8})
is generally different for different distributions of the pairs $\left(A_{ij},B_{ij}\right)$.
However, it always includes the coupling $Y$ in which the pairs $\left(A'_{ij},B'_{ij}\right)$
are stochastically independent across different $\left(i,j\right)$.
This coupling is referred to as an \emph{independent coupling}. Its
universal applicability leads to the common confusion of stochastic
unrelatedness with stochastic independence. But stochastic independence
is merely a special property of a joint distribution.

The non-uniqueness of the coupling (\ref{eq:coupling8}), rather than
being a hindrance, can be advantageously used in theoretical analysis.
According to the \emph{All-Possible-Couplings} principle formulated
in Refs. {[}\ref{Dzhafarov,-E.N.,-&LNCS13},\ref{Dzhafarov,-E.N.,-&Advances}{]}, 
\begin{quote}
a set of stochastically unrelated random variables is characterized
by the set of all possible couplings that can be imposed on them,
with no couplings being a priori privileged. 
\end{quote}
Thus, according to Ref. {[}\ref{Dzhafarov,-E.N.,-&2013}{]}, the set
of all possible couplings for (\ref{eq:8}) can be used to characterize
various constraints imposed on the joint distributions of $A_{ij},B_{ij}$
in (\ref{eq:all8}).

From the point of view of all possible couplings, the noncontextuality
assumption leading to the single-indexation of the spins, $A_{i},B_{j}$,
is equivalent to imposing an \emph{identity coupling} on the double-indexed
outputs in (\ref{eq:8}), i.e., creating a coupling (\ref{eq:coupling8})-(\ref{eq:coupling8 condition})
with the additional constraint{]} 
\begin{equation}
\begin{array}{cc}
\Pr\left[A_{11}^{\prime}=A_{12}^{\prime}\right]=1, & \Pr\left[A_{21}^{\prime}=A_{22}^{\prime}\right]=1,\\
\Pr\left[B_{11}^{\prime}=B_{21}^{\prime}\right]=1, & \Pr\left[B_{12}^{\prime}=B_{22}^{\prime}\right]=1.
\end{array}\label{eq:identity coupling}
\end{equation}
The Bell-type theorems {[}\ref{Bell,-J.:-On}-\ref{enu:Fine,-A.:-Hidden}{]}
tell us that this coupling exists if and only if both the no-signaling
condition is satisfied and the four observable joint distributions
of $\left(A_{ij},B_{ij}\right)$ satisfy the inequalities 
\begin{equation}
\begin{array}{c}
\left|\left\langle A_{11}B_{11}\right\rangle +\left\langle A_{12}B_{12}\right\rangle +\left\langle A_{21}B_{21}\right\rangle -\left\langle A_{22}B_{22}\right\rangle \right|\leq2,\\
\left|\left\langle A_{11}B_{11}\right\rangle +\left\langle A_{12}B_{12}\right\rangle -\left\langle A_{21}B_{21}\right\rangle +\left\langle A_{22}B_{22}\right\rangle \right|\leq2,\\
\left|\left\langle A_{11}B_{11}\right\rangle -\left\langle A_{12}B_{12}\right\rangle +\left\langle A_{21}B_{21}\right\rangle +\left\langle A_{22}B_{22}\right\rangle \right|\leq2,\\
\left|-\left\langle A_{11}B_{11}\right\rangle +\left\langle A_{12}B_{12}\right\rangle +\left\langle A_{21}B_{21}\right\rangle +\left\langle A_{22}B_{22}\right\rangle \right|\leq2,
\end{array}\label{eq:Bell}
\end{equation}
where $\left\langle \ldots\right\rangle $ denotes expected value.
Clearly, these inequalities do not have to be satisfied, and, in the
Alice-Bob paradigm, for some quadruples of settings $\left(\alpha_{1},\alpha_{2},\beta_{1},\beta_{2}\right)$,
these inequalities are contravened by quantum theory and experimental
data.

Therefore, we have to use double-indexing and consider couplings other
than the identity coupling (\ref{eq:identity coupling}). This is
the essence of the contextualization approach, when applied to the
Alice-Bob paradigm. In the conditionalization approach, discussed
next, one also uses what can be thought of as a version of double-indexation
(\emph{conditioning} on the two indices viewed as values of a random
variable), but instead of the couplings in the sense of (\ref{eq:coupling8})-(\ref{eq:coupling8 condition})
one uses a different theoretical construct, \emph{conditional couplings}.

\subsection*{Conditionalization and Conditional Couplings}

One of the simplest ways of creating stochastically unrelated random
variables is to consider a tree of possibilities, like this one: 
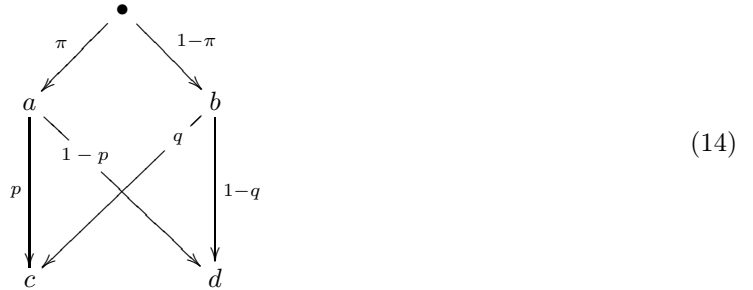
\begin{equation}
\begin{array}{c}
\xymatrix{ & \bullet\ar[dl]_{\pi}\ar[dr]^{1-\pi}\\
a\ar[dd]_{p}\ar[ddrr]|(0.3){\colorbox{white}{\scriptsize{\ensuremath{1-p}}}} &  & b\ar[dd]^{1-q}\ar[ddll]|(0.2){\colorbox{white}{\scriptsize{\ensuremath{q}}}}\\
\\
c &  & d
}
\end{array}\label{eq:tree}
\end{equation}
We have at the first stage outcomes $a$ and $b$, and according as
which of them is realized, the choice between $c$ and $d$ occurs
with generally different probabilities. We can consider $a$ and $b$
as two mutually exclusive conditions, and use them to label the two
random variables 
\begin{equation}
X_{a}=\begin{cases}
c & \textnormal{with probability }p,\\
d & \textnormal{with probability }1-p,
\end{cases}\quad X_{b}=\begin{cases}
c & \textnormal{with probability }q,\\
d & \textnormal{with probability }1-q.
\end{cases}
\end{equation}
Clearly, $X_{a}$ and $X_{b}$ here do not have a joint distribution:
e.g., no joint probability $\Pr\left[X_{a}=c,X_{b}=c\right]$ is defined
because there is no commonly acceptable meaning in which $X_{a}=c$
may ``co-occur'' with $X_{b}=c$. The two random variables here
are stochastically unrelated, in conformance with the Contextuality-by-Default
principle.

The All-Possible-Couplings principle leads us to consider all joint
distributions

\begin{equation}
\begin{array}{|c|c|c|}
\hline  & X'_{b}=c & X'_{b}=d\\
\hline X'_{a}=c & r & p-r\\
\hline X'_{a}=d & q-r & 1-p-q+r
\\\hline \end{array}\:,
\end{equation}
with 
\begin{equation}
\max\left(0,p+q-1\right)\leq r\leq\min\left(p,q\right).
\end{equation}
Each $r$ within this range defines a possible coupling $Y=\left(X'_{a},X'_{b}\right)$
for $X_{a}$ and $X_{b}$. In particular, the independent coupling,
with $r=pq$, is within the range, while the identity coupling, with
$\Pr\left[X_{a}=X_{b}\right]=1$, is possible if and only if $r=p=q$.

There is, however, a more traditional view of $X_{a}$ and $X_{b}$
in (\ref{eq:tree}). It consists in considering a joint distribution
of two random variables, $C$ and $X$, with the marginal distributions
\begin{equation}
C=\begin{cases}
a & \textnormal{with probability }\pi,\\
b & \textnormal{with probability }1-\pi,
\end{cases}\quad X=\begin{cases}
c & \textnormal{with probability }p\pi+q\left(1-\pi\right),\\
d & \textnormal{with probability }\left(1-p\right)\pi+\left(1-q\right)\left(1-\pi\right),
\end{cases}
\end{equation}
and with the joint distribution 
\begin{equation}
\begin{array}{|c|c|c|}
\hline  & X=c & X=d\\
\hline C=a & p\pi & \left(1-p\right)\pi\\
\hline C=b & q\left(1-\pi\right) & \left(1-q\right)\left(1-\pi\right)
\\\hline \end{array}\:.
\end{equation}
$X_{a}$ is then interpreted as $X$ given $C=a$, and analogously
for $X_{b}$. The conditional probabilities are computed as required,
\begin{equation}
p=\Pr\left[X=c\,\vert\, C=a\right],\quad q=\Pr\left[X=c\,\vert\, C=b\right].
\end{equation}

The idea suggested by this simple exercise is this: 
\begin{quote}
consider any set of stochastically unrelated random outputs labeled
by mutually exclusive conditions as if these conditions were values
of some random variable, and the outputs were values of another random
variable conditioned upon the values of the former. 
\end{quote}
We call this approach conditionalization. It may seem to provide a
simple alternative, within the framework of KPT, to considering all
couplings imposable on stochastically unrelated variables. We will
argue, however, that this alternative is not theoretically interesting.

Consider the conditionalization of our Alice-Bob paradigm. Denote,
for $i=1,2$ and $j=1,2$, 
\begin{equation}
\begin{array}{c}
p_{ij}=\Pr\left[A_{ij}=1,B_{ij}=1\right],\\
p_{i\cdot}=\Pr\left[A_{ij}=1\right],\\
p_{\cdot j}=\Pr\left[B_{ij}=1\right].
\end{array}
\end{equation}
Introduce a random variable $C$ with four values 
\[
c_{ij}=\left(\alpha_{i},\beta_{j}\right),\; i,j,\in\left\{ 1,2\right\} ,
\]
and a random variable $X=\left(A',B'\right)$ with four values 
\[
\left(1,1\right),\left(1,-1\right),\left(-1,1\right),\left(-1,-1\right).
\]
Form the tree of outcomes as shown below, using arbitrarily chosen
positive probabilities $\pi_{11},\pi_{12},\pi_{21},\pi_{22}$ (summing
to 1):

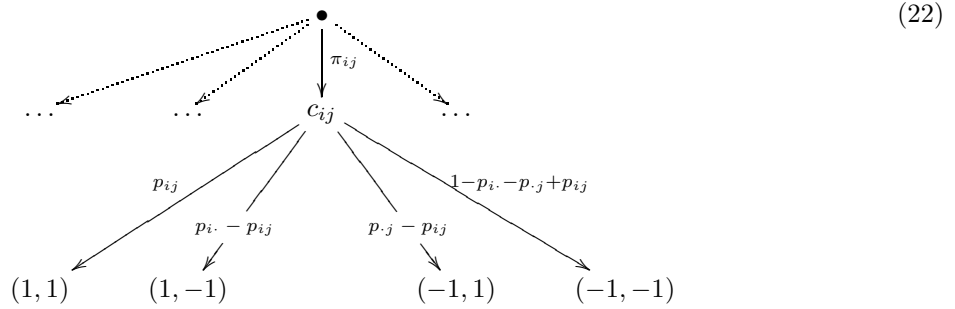
\begin{equation}
\xymatrix{ &  & \bullet\ar@{.>}[dll]\ar@{.>}[dl]\ar[d]^{\pi_{ij}}\ar@{.>}[dr]\\
\ldots & \ldots & c_{ij}\ar[ddll]_{p_{ij}}\ar[ddl]|(0.65){\colorbox{white}{\scriptsize{\ensuremath{p_{i\cdot}-p_{ij}}}}}\ar[ddr]|(0.65){\colorbox{white}{\scriptsize{\ensuremath{p_{\cdot j}-p_{ij}}}}}\ar[ddrr]^{1-p_{i\cdot}-p_{\cdot j}+p_{ij}} & \ldots\\
\\
\left(1,1\right) & \left(1,-1\right) &  & \left(-1,1\right) & \left(-1,-1\right)
}
\end{equation}
The conditionalization is completed by computing the joint distribution
of $C$ and $\left(A',B'\right)$: 
\begin{equation}
\begin{array}{|c|c|c|c|c|}
\hline  & \left(A',B'\right)=\left(1,1\right) & \left(1,-1\right) & \left(-1,1\right) & \left(-1,-1\right)\\
\hline \ldots & \ldots & \ldots & \ldots & \ldots\\
\hline C=\left(\alpha_{i},\beta_{j}\right) & p_{ij}\pi_{ij} & \left(p_{i\cdot}-p_{ij}\right)\pi_{ij} & \left(p_{\cdot j}-p_{ij}\right)\pi_{ij} & \left(1-p_{i\cdot}-p_{\cdot j}+p_{ij}\right)\pi_{ij}\\
\hline \ldots & \ldots & \ldots & \ldots & \ldots
\end{array}
\end{equation}
Clearly, we have constructed a random variable 
\begin{equation}
Z=\left(C,\left(A',B'\right)\right)\label{eq:simplest Z}
\end{equation}
such that 
\begin{equation}
\left(A',B'\right)\textnormal{ given }C=\left(\alpha_{i},\beta_{j}\right)\textnormal{ is distributed as }\left(A_{ij},B_{ij}\right).\label{eq:coupling3}
\end{equation}
This $Z$ can be called a conditional coupling for $\left(A_{ij},B_{ij}\right)$,
$i,j\in\left\{ 1,2\right\} $.

The conditionalization procedure does not have to claim the existence
of any ``true'' or unique distribution of $C$. One can freely concoct
this distribution, even if the conditions under which $A$ and $B$
are measured are chosen at will or according to a deterministic algorithm.

There are two interesting modifications of conditionalization, both
proposed in a recent paper by Avis, Fischer, Hilbert, and Khrennikov
{[}\emph{\ref{enu:D.-Avis,-P.}}{]}. Instead of the conditional coupling
$Z$ in (\ref{eq:simplest Z}), they consider 
\begin{equation}
Z'=\left(C,\left(A'_{1},A'_{2},B'_{1},B'_{2}\right)\right)
\end{equation}
such that
\begin{equation}
\left(A_{i}',B_{j}'\right)\textnormal{ given }C=\left(\alpha_{i},\beta_{j}\right)\textnormal{ is distributed as }\left(A_{ij},B_{ij}\right).
\end{equation}
In other words,
\begin{equation}
\Pr\left[A'_{i}=\pm1,B'_{j}=\pm1\,\vert\, C=\left(\alpha_{i},\beta_{j}\right)\right]=\Pr\left[A_{ij}=\pm1,B_{ij}=\pm1\right].\label{eq:Avis0}
\end{equation}
This does not yet define the conditional probabilities for all possible
values of $\left(A'_{1},A'_{2},B'_{1},B'_{2}\right)$. Avis et al.
describe two ways of defining them. 

In one of them $A'_{1},A'_{2},B'_{1},B'_{2}$ have two possible values
each, $\pm1$, and
\begin{equation}
\Pr\left[A'_{i}=a,B'_{j}=b,\, A'_{3-i}=a',B'_{3-j}=b',|\, C=\left(\alpha_{i},\beta_{j}\right)\right]=\frac{1}{4}\Pr\left[A_{ij}=a,B_{ij}=b\right].\label{eq: conditional v1}
\end{equation}
That is, the probability of $\left(A'_{i}=a,B'_{j}=b\right)$ at $C=\left(\alpha_{i},\beta_{j}\right)$
is evenly partitioned among the four values of the ``irrelevant''
pair $\left(A'_{3-i},B'_{3-j}\right)$. It is easy to see that one
could as well use any other partitioning: 
\begin{equation}
\Pr\left[A'_{i}=a,B'_{j}=b,\, A'_{3-i}=a',B'_{3-j}=b',|\, C=\left(\alpha_{i},\beta_{j}\right)\right]=t_{ij}\left(a',b'\right)\Pr\left[A_{ij}=a,B_{ij}=b\right],
\end{equation}
with nonnegative $t_{ij}\left(a',b'\right)$ subject to 
\[
t_{ij}\left(1,1\right)+t_{ij}\left(1,-1\right)+t_{ij}\left(-1,1\right)+t_{ij}\left(-1,-1\right)=1,\; i,j\in\left\{ 1,2\right\} .
\]
Now, for any distribution of $C$ with non-zero values of $\Pr\left[C=\left(\alpha_{i},\beta_{j}\right)\right]$,
the joint distribution of $\left(C,A'_{1},A'_{2},B'_{1},B'_{2}\right)$
is well-defined. 

Another way of implementing (\ref{eq:Avis0}) described in Ref. {[}\emph{\ref{enu:D.-Avis,-P.}}{]}
is to allow each of $A'_{1},A'_{2},B'_{1},B'_{2}$ to attain a third
value, say, $0$, in addition to $\pm1$. This third value can be
interpreted as ``is not defined.'' It is postulated then that

\begin{equation}
\Pr\left[A'_{i}=a,B'_{j}=b,\, A'_{3-i}=a',B'_{3-j}=b',|\, C=\left(\alpha_{i},\beta_{j}\right)\right]=\begin{cases}
\Pr\left[A_{ij}=a,B_{ij}=b\right] & \textnormal{if }a\not=0,b\not=0,a'=b'=0,\\
0 & \textnormal{otherwise}.
\end{cases}\label{eq: conditional v2}
\end{equation}
Again, it is easy to see that the joint distribution of $\left(C,A'_{1},A'_{2},B'_{1},B'_{2}\right)$
is well-defined and satisfies (\ref{eq:Avis0}) for any distribution
of $C$ with non-zero values of $\Pr\left[C=\left(\alpha_{i},\beta_{j}\right)\right]$.

\section*{Discussion: Comparing the Two Approaches}

Conditionalization and contextualization achieve the same goal ---
``sewing together'' stochastically unrelated random variables within
the confines of KPT. But the similarity ends here. Consider, e.g.,
the Alice-Bob experiment in which both Alice and Bob use some random
generators to choose between two possible measurement directions.
Clearly then $C$ is objectively a random variable, and a joint distribution
of $\left(A,B\right)$ and $C$ objectively exists. Put differently,
in this case $\left(A',B'\right)$ given $C=\left(\alpha_{i},\beta_{j}\right)$
in (\ref{eq:coupling3}) is simply equal to $\left(A_{ij},B_{ij}\right)$.

However, whether $C$ is objectively a random variable or a distribution
for the settings is invented, the quantum-mechanical analysis of the
situation begins with computing the (conditional) distributions of
$\left(A,B\right)$ at different settings. The distribution of $C$
in no way advances our understanding of how $\left(A_{ij},B_{ij}\right)$
for different $\left(i,j\right)$ are related to each other.

Thus, we know that the entangled spin-$\nicefrac{1}{2}$ particles
are subject to Tsirelson's inequalities {[}\ref{Cirel'son,-B.S.:-Quantum}{]}
\begin{equation}
\begin{array}{c}
\left|\left\langle A_{11}B_{11}\right\rangle +\left\langle A_{12}B_{12}\right\rangle +\left\langle A_{21}B_{21}\right\rangle -\left\langle A_{22}B_{22}\right\rangle \right|\leq2\sqrt{2},\\
\left|\left\langle A_{11}B_{11}\right\rangle +\left\langle A_{12}B_{12}\right\rangle -\left\langle A_{21}B_{21}\right\rangle +\left\langle A_{22}B_{22}\right\rangle \right|\leq2\sqrt{2},\\
\left|\left\langle A_{11}B_{11}\right\rangle -\left\langle A_{12}B_{12}\right\rangle +\left\langle A_{21}B_{21}\right\rangle +\left\langle A_{22}B_{22}\right\rangle \right|\leq2\sqrt{2},\\
\left|-\left\langle A_{11}B_{11}\right\rangle +\left\langle A_{12}B_{12}\right\rangle +\left\langle A_{21}B_{21}\right\rangle +\left\langle A_{22}B_{22}\right\rangle \right|\leq2\sqrt{2}.
\end{array}\label{eq:Tsirelson}
\end{equation}
We also know that if the two particles were not entangled, they would
be subject to the Bell-CH-Fine inequalities (\ref{eq:Bell}). The
difference between these two constraints is not reflected in the ``true''
distribution of $C$, if it exists, nor is it implied by or can in
any way restrict the possible choices of ``imaginary'' distributions
of $C$. In fact, the only restriction imposed on the distribution
of $C$, a universal one, is that none of the conditions should have
probability zero, because this would make the conditional probabilities
undefined. Moreover, the set of possible conditional couplings is
the same whether the no-signaling condition is or is not satisfied.

Although in this discussion we assumed that conditionalization was
implemented in its simplest version, (\ref{eq:simplest Z})-(\ref{eq:coupling3}),
our arguments and conclusions apply verbatim to the modifications
proposed in Ref. {[}\emph{\ref{enu:D.-Avis,-P.}}{]} and described
at the end of the previous section. The conditional distributions
of $A'_{1},A'_{2},B'_{1},B'_{2}$ for the four values of $C$ in (\ref{eq: conditional v1})
and (\ref{eq: conditional v2}) are uniquely determined by the observed
distributions of the four pairs $\left(A_{ij},B_{ij}\right)$. But
whatever these distributions, they can be paired with any distribution
of $C$, provided none of its values has zero probability. 

All of this stands in a clear contrast to the analysis of all possible
couplings (\ref{eq:coupling8}) in the contextualization approach
{[}\ref{Dzhafarov,-E.N.,-&2013}-\ref{Dzhafarov,-E.N.,-&New}{]}.
In this approach we can ask various questions about the compatibility
of couplings with various constraints known to hold for the observable
joint distributions. Thus, we may ask about the \emph{fitting set}
of couplings for a given constraint (say, Bell or Tsirelson inequalities),
i.e., the couplings that are compatible with the spin distributions
subject to the constraint. We can also ask about the \emph{forcing
set} of couplings, those compatible only with the spin distributions
subject to a given constraint. Or we can conjoin the two questions
and ask about the \emph{equivalent set} of couplings, those compatible
with and only with the spin distributions subject to the constraint.
The answers to such questions will be different for different constraints
being considered.

Since the four observed joint distributions of $\left(A'_{ij},B'_{ij}\right)$
in (\ref{eq:coupling8 condition}) are themselves part of the couplings
(\ref{eq:coupling8}), the questions above are only interesting if
they are formulated in terms of the unobservable parts of the couplings.
In the examples below we characterize the couplings in terms of the
\emph{connections} {[}\ref{Dzhafarov,-E.N.,-&2013},\ref{Dzhafarov,-E.N.,-&LNCS13},\ref{Dzhafarov,-E.N.,-&New}{]},
which are the (unobservable) pairs 
\begin{equation}
\left(A'_{11},A'_{12}\right),\left(A'_{21},A'_{22}\right),\left(B'_{11},B'_{21}\right),\left(B'_{12},B'_{22}\right).\label{eq:connections}
\end{equation}
The diagram below shows the connections in their relation to the pairs
whose joint distributions are known from observations (compare with
diagram (\ref{eq:all8})): 
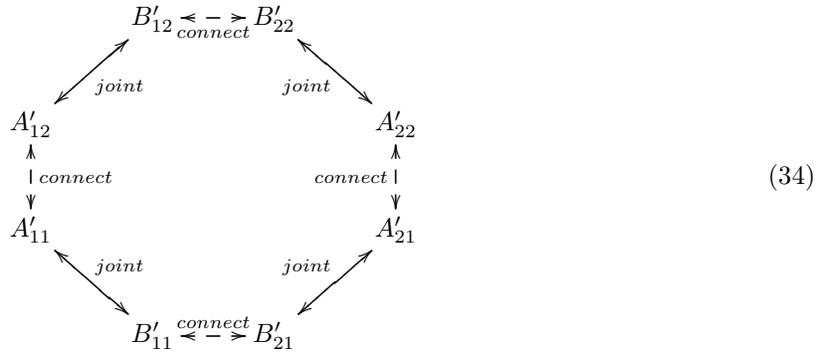
\begin{equation}
\begin{array}{c}
\xymatrix{ & B'_{12}\ar[dl]\ar@{-->}[r]_{connect} & B'_{22}\ar[dr]_{joint}\ar@{-->}[l]\\
A'_{12}\ar[ur]_{joint}\ar@{-->}[d]^{connect} &  &  & A'_{22}\ar[ul]\ar@{-->}[d]_{connect}\\
A'_{11}\ar[dr]^{joint}\ar@{-->}[u] &  &  & A'_{21}\ar[dl]\ar@{-->}[u]\\
 & B'_{11}\ar[lu]\ar@{-->}[r] & B'_{21}\ar[ur]^{joint}\ar@{-->}[l]_{connect}
}
\end{array}
\end{equation}

Let us assume that the probability of spin-up ($+1$) outcome for
every (spin-$\nicefrac{1}{2}$) particle in the Alice-Bob paradigm
is $\nicefrac{1}{2}$. (As shown in Ref. {[}\ref{enu:Masanes,-Ll.,-Acin,}{]},
this can always be achieved by a simple procedural modification of
the canonical Alice-Bob experiment.) This assumption is, of course,
in compliance with the no-signaling condition, which therefore can
be omitted from all formulations below.

We know {[}\ref{Dzhafarov,-E.N.,-&2013}{]} that the following two
statements about connections are equivalent:
\begin{quote}
($S_{1}$) a vector of connections (\ref{eq:connections}) is compatible
\emph{with and only with} those distributions of $\left(A_{ij},B_{ij}\right)$,
$i,j\in\left\{ 1,2\right\} $, that satisfy the Bell-CH-Fine inequalities
(\ref{eq:Bell}); 
\[
\textnormal{is equivalent to}
\]

($S_{1}'$) a vector of connections (\ref{eq:connections}) is such
that
\begin{equation}
\left\langle A_{11}A_{12}\right\rangle =\pm1,\left\langle A_{21}A_{22}\right\rangle =\pm1,\left\langle B_{11}B_{21}\right\rangle =\pm1,\left\langle B_{21}B_{22}\right\rangle =\pm1,
\end{equation}
where the number of + signs among the four expected values is 4,2,
or 0. 
\end{quote}
The equivalence of these two statements is an expanded version of
Fine's theorem {[}\ref{enu:Fine,-A.:-Hidden}{]}, whose formulation
in the language of connections is: the identity connections, those
with 
\begin{equation}
\left\langle A_{11}A_{12}\right\rangle =\left\langle A_{21}A_{22}\right\rangle =\left\langle B_{11}B_{21}\right\rangle =\left\langle B_{21}B_{22}\right\rangle =1,
\end{equation}
are only compatible with distributions of $\left(A_{ij},B_{ij}\right)$
satisfying the Bell-CH-Fine inequalities; and if these inequalities
hold, then $\left(A_{ij},B_{ij}\right)$ can be coupled by means of
the identity connections.

We also know {[}\ref{Dzhafarov,-E.N.,-&2013}{]} that the following
two statements about connections are equivalent: 
\begin{quote}
($S_{2}$) a vector of connections (\ref{eq:connections}) is compatible
\emph{with and only with} those distributions of $\left(A_{ij},B_{ij}\right)$,
$i,j,\in\left\{ 1,2\right\} $, that satisfy the Tsirelson inequalities
(\ref{eq:Tsirelson}); 
\[
\textnormal{is equivalent to}
\]

($S'_{2}$) a vector of connections (\ref{eq:connections}) is such
that
\begin{equation}
\max\left\{ \pm\left\langle A_{11}A_{12}\right\rangle \pm\left\langle A_{21}A_{22}\right\rangle \pm\left\langle B_{11}B_{21}\right\rangle \pm\left\langle B_{21}B_{22}\right\rangle :\textnormal{number of +'s is even}\right\} =2\left(3-\sqrt{2}\right)
\end{equation}
and 
\begin{equation}
\max\left\{ \pm\left\langle A_{11}A_{12}\right\rangle \pm\left\langle A_{21}A_{22}\right\rangle \pm\left\langle B_{11}B_{21}\right\rangle \pm\left\langle B_{21}B_{22}\right\rangle :\textnormal{number of +'s is odd}\right\} \leq2.
\end{equation}

\end{quote}
We see that although the expectations $\left\langle A_{i1}A_{i2}\right\rangle $
and $\left\langle B_{1j}B_{2j}\right\rangle $ for the connections
are not observable, they provide a theoretically meaningful way of
characterizing the way in which the stochastically unrelated and observable
$\left(A_{ij},B_{ij}\right)$ are being ``sewn together.'' And these
ways are different for the Bell-CH-Fine and Tsirelson inequalities.

What can contextualization tell us about the basic predictions of
the quantum theory for the Alice-Bob experiment? The theory tells
us that, for $i=1,2$ and $j=1,2$, 
\begin{equation}
\left\langle A_{ij}B_{ij}\right\rangle =-\langle\alpha_{i}\,\vert\,\beta_{j}\rangle,\label{eq:quantum correaltions}
\end{equation}
where $\langle\alpha_{i}|\beta_{j}\rangle$ is the dot product of
two unit vectors. It can be shown {[}\ref{Landau,-L.J.-(1988).}-\ref{Kujala,-J.V.,-&}{]}
that the four expectations $\left\langle A_{ij}B_{ij}\right\rangle $
can be presented in the form (\ref{eq:quantum correaltions}) using
a quadruple of setting $\left(\alpha_{1},\alpha_{2},\beta_{1},\beta_{2}\right)$
if and only if 
\begin{equation}
\begin{array}{c}
\left|\arcsin\left\langle A_{11}B_{11}\right\rangle +\arcsin\left\langle A_{12}B_{12}\right\rangle +\arcsin\left\langle A_{21}B_{21}\right\rangle -\arcsin\left\langle A_{22}B_{22}\right\rangle \right|\leq\pi,\\
\left|\arcsin\left\langle A_{11}B_{11}\right\rangle +\arcsin\left\langle A_{12}B_{12}\right\rangle -\arcsin\left\langle A_{21}B_{21}\right\rangle +\arcsin\left\langle A_{22}B_{22}\right\rangle \right|\leq\pi,\\
\left|\arcsin\left\langle A_{11}B_{11}\right\rangle -\arcsin\left\langle A_{12}B_{12}\right\rangle +\arcsin\left\langle A_{21}B_{21}\right\rangle +\arcsin\left\langle A_{22}B_{22}\right\rangle \right|\leq\pi,\\
\left|-\arcsin\left\langle A_{11}B_{11}\right\rangle +\arcsin\left\langle A_{12}B_{12}\right\rangle +\arcsin\left\langle A_{21}B_{21}\right\rangle +\arcsin\left\langle A_{22}B_{22}\right\rangle \right|\leq\pi.
\end{array}\label{eq:quantum}
\end{equation}
These inequalities are ``sandwiched'' between the Bell-CH-Fine ones
and Tsirelson ones. That is, they are implied by the former and imply
the latter. It is shown in Ref. {[}\ref{Dzhafarov,-E.N.,-&New}{]}
that 
\begin{quote}
($S_{3}$) there is \emph{no} vector of connections (\ref{eq:connections})
that is compatible \emph{with and only with} those distributions of
$\left(A_{ij},B_{ij}\right)$, $i,j,\in\left\{ 1,2\right\} $, that
satisfy the quantum inequalities (\ref{eq:quantum}). 
\end{quote}
Moreover, this negative statement still holds if one replaces the
connections (\ref{eq:connections}) with any other subsets of (\ref{eq:coupling8}),
e.g., 
\begin{equation}
\left(A_{11}^{\prime},A_{12}^{\prime},A_{21}^{\prime},A_{22}^{\prime}\right),\left(B_{11}^{\prime},B_{12}^{\prime},B_{21}^{\prime},B_{22}^{\prime}\right).
\end{equation}
No distributions of such subsets are compatible with and only with
those distributions of $\left(A_{ij},B_{ij}\right)$ that satisfy
the quantum inequalities (\ref{eq:quantum}).

The investigation of the forcing set of couplings provides additional
insights into the special nature of quantum mechanics. The result
we have {[}\ref{Dzhafarov,-E.N.,-&New}{]} says that the following
two statements about connections are equivalent (note the change from
``with and only with'' of the previous statements to ``only with''):
\begin{quote}
($S_{4}$) a vector of connections (\ref{eq:connections}) is compatible
\emph{only with} those distributions of $\left(A_{ij},B_{ij}\right)$,
$i,j\in\left\{ 1,2\right\} $, that satisfy the quantum inequalities
(\ref{eq:quantum});
\[
\textnormal{is equivalent to}
\]

($S'_{4}$) a vector of connections (\ref{eq:connections}) is compatible
\emph{only with} those distributions of $\left(A_{ij},B_{ij}\right)$,
$i,j\in\left\{ 1,2\right\} $, that satisfy the Bell-CH-Fine inequalities
(\ref{eq:Bell}). 
\end{quote}
In other words, a choice of connections can force all $\left(A_{ij},B_{ij}\right)$
compatible with them to comply with quantum mechanics only in the
form of their compliance with classical mechanics.

\section*{Conclusion}

The examples just given should suffice to illustrate the point made:
while both contextualization and conditionalization embed any input-output
relation into the framework of KPT, only contextualization provides
a useful tool for understanding the nature of various constraints
imposed on the observable joint distributions (one could say also,
for different types and levels of contextuality). Conditionalization
is uninformative, as any distribution of the conditions is compatible
with any distributions of the conditional random variables.

\end{document}